\documentclass[12pt]{amsart}

\title{Minimal models of $L^*$}

\author{Ying Zong}

\address{}
\email{zongying.upenn@gmail.com}

\date{} 

\begin{document}

\maketitle

\section{Introduction}

\smallskip

Let $S$ be an algebraic space, $A$ an $S$-abelian algebraic space, and $L$ an invertible module on $A$. We suppose that $L$ is fiberwise numerically trivial. This assumption, according to the theorem of square, is equivalent to that, on the $\mathbf{G}_m$-torsor $L^*:=\underline{\mathrm{Isom}}_A(\mathcal{O}_A, L)$, there exists, locally for the \'{e}tale topology of $S$, a compatible structure of extension of $A$ by $\mathbf{G}_m$:
\[1\to \mathbf{G}_m\to L^*\to A\to 1\]

\smallskip

{\bf Definition.} --- \emph{A minimal model of $L^*\times \mathbf{Z}/n\mathbf{Z}$, $n\geq 1$, consists of a pair $(f, \varphi)$, where }:
\smallskip

\emph{a) $f: X\to S$ is a proper flat finitely presented $S$-algebraic space, with geometrically connected fibers, such that, for each geometric point $s$ of $S$, there is an isomorphism
\[(f\times_Ss)^{!}k(s)\simeq \mathcal{O}[d]\] in $D(f^{-1}(s), \mathcal{O})$, $d=\mathrm{dim}\ f^{-1}(s)$.}

\smallskip

\emph{b) $\varphi$ is an $S$-pseudo-isomorphism of $X$ to $L^*\times\mathbf{Z}/n\mathbf{Z}$, which is defined and \'{e}tale at all points $x$ of codimension $\leq 1$ in $f^{-1}(f(x))$ where $f$ is smooth.}

\smallskip

\smallskip

This notion, relevant in the study of fibers of lagrangian fibrations on symplectic manifolds, \emph{cf}. \cite{weierstrass} Prop. 3, is motivated by \emph{N\'{e}ron polygons} of Deligne-Rapoport \cite{deligne2}, the starting point for investigating boundaries of various moduli of elliptic curves. My purpose, more birational, was to see whether, even under the weakest possible assumption on the pseudo-isomorphism $\varphi$, the following example, which we call \emph{standard} $n$-\emph{gon} analogous to \emph{loc.cit.}, exhausts in fact all minimal models, at least fiberwise.

\smallskip

\smallskip

{\bf Example.} --- Let the projective line over $\mathrm{Spec}(\mathbf{Z})$ be $P$, let it be pinched with the coordinate sections $T=0, \infty$ identified, and let $c: C\to \mathrm{Spec}(\mathbf{Z})$ be the resulting nodal curve. One can make explicit, via the normalization $P\to C$, the dualizing module $H^{-1}c^{!}\mathcal{O}$; it has as $\mathcal{O}_C$-basis the logarithmic differential $d_{P/\mathbf{Z}}(T)/T$, which is of residue $1$ and $-1$ along $T=0$ and $T=\infty$ respectively. And one calculates the group $H^1(C, \mathbf{Z})$ of the isomorphism classes of $\mathbf{Z}$-torsors, in the sense of \'{e}tale topology, on $C$; it is infinite cyclic, with a generator $C_{\infty}$ represented by an unending chain of projective lines, any two lines in succession glued along the section $\infty$ on the former and along the section $0$ on the latter. The multiplicative group $\mathbf{G}_m$ acts canonically, compatibly on $C$ and $C_{\infty}$ as well as all intermediate quotients $C_{\infty}/n\mathbf{Z}$. The compactification
\[C_{nL}:=L^*\stackrel{\mathbf{G}_m}{\wedge}C_{\infty}/n\mathbf{Z}\] of $L^*\times\mathbf{Z}/n\mathbf{Z}$ is a minimal model.

\smallskip

\smallskip

It turns out, when $A$ is non-trivial, models other than the standard do exist. We call them $n$-\emph{chains}, all of which, fiber by fiber, are obtained from $C_{nL}$ by normalization and gluing and pinching (Artin \cite{artin} 6.1).

\smallskip

\smallskip

{\bf Theorem.} --- \emph{Let $(f, \varphi)$ be a minimal model of $L^*\times\mathbf{Z}/n\mathbf{Z}$. Then }:
\smallskip

\emph{i) The $S$-pseudo-isomorphism $\varphi^{-1}$ is everywhere defined and an isomorphism onto the domain of definition of $\varphi$. The geometric number of connected components of $f^{-1}(s)-\mathrm{Im}(\varphi_s^{-1})$, $s\in S$, is $n$ or $n+1$, and locally constant on $S[1/2]$.}

\smallskip

\emph{ii) Assume $S=\mathrm{Spec}(k)$, $k$ an algebraically closed field, and $X-\mathrm{Im}(\varphi^{-1})$ of $n$ connected components. Then $X$ is $C_{nL}$.}

\smallskip

\emph{iii) Assume $S=\mathrm{Spec}(k)$, $k$ an algebraically closed field, and $X-\mathrm{Im}(\varphi^{-1})$ of $n+1$ connected components. Let 
\[P_L:=L^*\stackrel{\mathbf{G}_m}{\wedge}\mathbf{P}^1\] and $\{0, \infty\}$ its distinguished pair of sections over $A$, and let $Q_{nL}$ be the chain obtained from $P_L\times\mathbf{Z}/n\mathbf{Z}$ by gluing the $i$-th section $\infty_i$ with the $(i+1)$-th section $0_{i+1}$, $i=1,\cdots, n-1$. Then $X$ is pinched from $Q_{nL}$ along finite $k$-morphisms 
\[\rho^1: A=0_1\to \Sigma_1,\ \rho^n: A=\infty_n\to \Sigma_n\] both of generic degree $2$, with $\Sigma_i$ proper over $\mathrm{Spec}(k)$, normal, integral, Cohen-Macaulay, of dualizing module $K_{\Sigma_i}=\rho^i_{*}\mathcal{O}_A/\mathcal{O}_{\Sigma_i}$, $i=1, n$. Fix either $\rho: A\to\Sigma$ and a volume form $\mu$ on $A/k$.}

\smallskip

\emph{When $k(A)/k(\Sigma)$ is inseparable, the relative Frobenius of $A/k$ factors through $\rho$ as
\[F: A\stackrel{\rho}{\to}\Sigma\stackrel{\varepsilon}{\to}A'\] with $\varepsilon$ finite flat, and $\mathrm{Tr}_{\rho}(\mu)=0$ in $H^0(\Sigma, K_{\Sigma})$.}

\smallskip

\emph{When $k(A)/k(\Sigma)$ is separable, $\Sigma$ is the quotient $A/\sigma$ by an involutive automorphism $\sigma$ of $A$, and $\sigma^*\mu=-\mu$. In characteristic $2$, the fixed point set of $\sigma$ is either empty or purely of codimension $1$ in $A$.}

\smallskip

\smallskip

The essential of the proof is arguing that, over the spectrum of an algebraically closed field, $X$ shares normalization with $C_{nL}$, see I) \emph{p}) below. The fact that $\varphi^{-1}$ is everywhere defined and an isomorphism onto $\mathrm{Dom}(\varphi)$ is proved in III) \emph{d}). That the geometric number of connected components of $f^{-1}(s)-\mathrm{Im}(\varphi^{-1}_s)$, $s\in S[1/2]$, is locally constant is proved in III) \emph{h}). 

\smallskip

\section{Proof of the theorem}

\smallskip

I) \emph{Assume first $n=1$ and that $S$ is the spectrum of an algebraically closed field $k$.}

\smallskip

Let $d=\mathrm{dim}\ X$. Since $X$ is pseudo-isomorphic to $L^*$, $X$ is integral, and $A$ is of dimension $d-1$. Let $t$ denote the generic point of $A$.

\smallskip

Let $U$ consist of all points of $X$ where $X$ is normal. Let $\rho: X'\to X$ be the normalization of $X$, and let $U$ be via $\rho$ identified with $\rho^{-1}(U)$.
\smallskip 

Write $I=H^0\rho^{!}\mathcal{O}$ for the conductor of $X'/X$, coherent ideal in $\mathcal{O}_{X'}$ and in $\mathcal{O}_X$, by which one defines the preferred closed sub-algebraic space structure of $Z'=X'-U$ (resp. $Z=X-U$) in $X'$ (resp. $X$).
\smallskip

The hypothesis that $f^{!}k\simeq \mathcal{O}[d]$ in $D(X, \mathcal{O})$ says in particular that $X$ is Cohen-Macaulay. As $I\simeq H^{-d}(f\rho)^{!}k$, $I$ is dualizing at precisely those points of $X'$ where $X'$ is Cohen-Macaulay.

\smallskip

Let $P$ be the projective line over $\mathrm{Spec}(k)$, and
\[P_L:=L^*\stackrel{\mathbf{G}_m}{\wedge}P\] and $\varphi'$ the pseudo-isomorphism of $X'$ to $P_L$ induced by $\varphi$. 
\smallskip

Let $D'$ consist of all points $x'$ of $X'$ where $X'$ is $k$-smooth and $\varphi'$ is defined. One has $\mathrm{codim}(X'-D', X')\geq 2$.
\smallskip

Write $C_L$ for $C_{1L}$, and $\pi: C_L\to\mathrm{Spec}(k)$ its structural morphism. The normalization of $C_L$ is $P_L$. Let $J$ be the conductor of $P_L/C_L$. In $P_L$, $V(J)$ is formed of two disjoint $A$-sections $0, \infty$. 
\smallskip

A basis of $H^{-d}\pi^{!}k$ can be thought of as a non-zero global $d$-form on $P_L/\mathrm{Spec}(k)$ with simple poles along $V(J)$. Such a $d$-form $\lambda$ has residue along either section $i=0, \infty$ a volume form $\mathrm{Res}_i(\lambda)\in H^0(A, \Omega^{d-1}_{A/k})$, and
\[\mathrm{Res}_0(\lambda)+\mathrm{Res}_{\infty}(\lambda)=0\]

\smallskip

\emph{a) The pseudo-morphism $\varphi|U$ is defined on all of $U$ }:

\smallskip

By hypothesis of minimal model, $\varphi$ is defined at all codimension $\leq 1$ points of $U$. Now apply Weil \cite{neron} 4.4/1, $L^*$ being a smooth $k$-algebraic group.

\smallskip

\emph{b) The morphism $\varphi|U$ is \'{e}tale }:

\smallskip

It is \'{e}tale at every codimension $\leq 1$ point of $U$ by hypothesis, hence \'{e}tale throughout by purity of branch locus \cite{almost} 2.4.

\smallskip

\emph{c) The birational morphism $\varphi|U$ induces an open immersion of $U$ into $L^*$ }: 

\smallskip

This follows by \emph{b)} and Zariski's Main Theorem. By \emph{c)} we identify $U$ with its image in $L^*$. 

\smallskip

\smallskip

\emph{d) The pseudo-morphism of $X'$ to $A$, composition of $\rho$, of $\varphi$ and of the projection $L^*\to A$, is everywhere defined }:

\smallskip

This follows by Weil \cite{neron} 4.4/1, since $X'$ is normal and $A$ an abelian variety.

\smallskip

\emph{e) The fiber of the morphism $X'\to A$ in d) over the generic point $t$ of $A$ is a projective $t$-line, on which $Z'_t$ consists of two distinct $t$-rational points }:

\smallskip

The fiber $X'_t$ is proper over $t$, normal, in which $U_t$ is open dense. The open immersion
\[\varphi': X'_t-Z'_t=U_t\hookrightarrow L^*_t=P_{Lt}-\{0, \infty\}\] uniquely extends to an isomorphism of $X'_t$ with $P_{Lt}\simeq \mathbf{P}^1_t$, under which the image of $Z'_t$ contains $\{0, \infty\}$. So $X'_t$ is a projective $t$-line, $I_t$ \emph{the} dualizing module, and $Z'_t$ an anti-canonical divisor, on $X'_t$.  Being of degree $2$ over $t$, $Z'_t$ consists of two distinct $t$-rational points.

\smallskip

\smallskip

Write $Z'_t=\{0, \infty\}$, and let $\Sigma'_i$ be the closed image of $i$, $i=0, \infty$, in $Z'$.

\smallskip

Note that the formalism of duality implies 
\[I|D'\simeq (\varphi'|D')^{!}J\] which, in terms of $1$-codimensional cycles, amounts to a linear equivalence 
\[-Z'|D'\sim -\Sigma'_0|D'-\Sigma'_{\infty}|D'+E'\] where $E'$ is $(\varphi'|D')$-exceptional.

\smallskip

\smallskip

\emph{f) One has $E'\geq 0$ }:

\smallskip

For, $P_L$ is smooth over $\mathrm{Spec}(k)$, \emph{cf}. the proof of \cite{almost} 4.5.

\smallskip

\smallskip

\emph{g) One has
$Z'|D'=(\Sigma'_0+\Sigma'_{\infty})|D'$, $E'=0$ }:
\smallskip

The cycle
\[Z'|D'-\Sigma'_0|D'-\Sigma'_{\infty}|D'+E'\] is $\geq 0$ and linearly equivalent to $0$, therefore is $0$, since 
\[H^0(D', \mathcal{O})=H^0(X', \mathcal{O})=k\]

\smallskip

\emph{h) The pseudo-isomorphism of $X'-(\Sigma'_0\cup\Sigma'_{\infty})$ to $L^*$ induced by $\varphi'$ is everywhere defined and an open immersion. In particular, the closed immersion $\Sigma'_0\cup\Sigma'_{\infty}\hookrightarrow Z'$ is surjective. The restriction 
\[\rho: \Sigma'_0\cup\Sigma'_{\infty}\to Z\] is surjective, so $Z$ either is irreducible or consists of two irreducible components }:

\smallskip

By \emph{g}) the complement of $U=X'-Z'$ in $X'-(\Sigma'_0\cup\Sigma'_{\infty})$ is of codimension $\geq 2$. Since $U\hookrightarrow L^*$ is an open immersion, the claim is immediate by the combination of Weil \cite{neron} 4.4/1, purity of branch locus \cite{almost} 2.4 and Zariski's Main Theorem.

\smallskip

\smallskip

\emph{i) The birational morphism $\varphi'|D'$ induces an open immersion of $D'$ into $P_L$ }:

\smallskip

By Zariski's Main theorem and by purity of branch locus \cite{almost} 2.4, it suffices to show that $\varphi'|D'$ is \'{e}tale at all codimension $\leq 1$ points of $D'$. According to \emph{f}), the $(\varphi'|D')$-exceptional cycle $E'$ is $0$, hence \emph{i}). 

\smallskip

\smallskip

\emph{j) The domain of definition of $\varphi'$ is $D'$ }:

\smallskip

By definition, $D'$ consists of all points of $\mathrm{Dom}(\varphi')$ where $X'$ is $k$-smooth. As $X'-D'$ is of codimension $\geq 2$ in $X'$, and $\varphi'|D'$ an open immersion \emph{i}), it follows by purity of branch locus \cite{almost} 2.4 that $\varphi'$ is \'{e}tale on $\mathrm{Dom}(\varphi')$. \emph{A priori}, $\mathrm{Dom}(\varphi')$ is $k$-smooth, hence \emph{j}).

\smallskip

Let $D'$ be via $\varphi'|D'$ identified as an open sub-scheme of $P_L$ so that $I|D'=J|D'$.

\smallskip

\emph{k) The fiber of $X'\to A$ over every codimension $\leq 1$ point of $A$ is a projective line }:

\smallskip

This has been shown over the generic point $t$ of $A$. Assume $s$ is a codimension $1$ point of $A$. Clearly, $X'\to A$ is flat around $X'_s$. Since $\mathrm{codim}(X'-D', X')\geq 2$, $D'_s$ is schematically dense in $X'_s$. So $X'_s$ is geometrically integral and, as a flat specialization of $X'_t$, is a projective $s$-line itself. 

\smallskip

\smallskip

\emph{l) Let $V'$ consist of all points of $X'$ where $X'$ is smooth over $\mathrm{Spec}(k)$. Then $\mathrm{codim}(X'-V', X')\geq 3$ }:

\smallskip

Otherwise, there would exist a maximal point $x'$ of $X'-V'$, of codimension $2$ in $X'$, with image $s$ of codimension $>1$ in $A$ by \emph{k}). 
\smallskip

Write $\Gamma$ (resp. $\overline{\Gamma}$) for the graph of $\varphi'$ (resp. the closed image of $\Gamma$ in $X'\times_AP_L$). Let 
\[\gamma: \widetilde{\Gamma}\to\overline{\Gamma}\] be a proper morphism from a normal integral scheme $\widetilde{\Gamma}$ such that $\gamma$ is an isomorphism over $D'\subset X'$ and such that $\widetilde{\Gamma}$ is regular above $x'$. Such partial resolution exists since excellent two dimensional singularities are canonically resolvable. 
\smallskip

Call $h$ (resp. $\psi$) the projection $\widetilde{\Gamma}\to\overline{\Gamma}\to X'$ (resp. $\widetilde{\Gamma}\to\overline{\Gamma}\to P_L$). 
\smallskip

The ideal $I$ being dualizing at $x'$, let $h^{!}I$ and $\psi^{!}J$ be along $h^{-1}(x')$ identified via the unique isomorphism that is compatible with $I|D'=J|D'$. 
\smallskip

Each maximal point $v$ of $h^{-1}(x')$ is of codimension $1$ in $\widetilde{\Gamma}$ and, as $s$ is of codimension $>1$ in $A$, is $\psi$-exceptional. One has 
\[(\psi^{!}J)_v=\mathcal{O}_v(r_v v)\] for an integer $r_v\geq 0$, since $P_L$ is regular (\emph{cf}. the proof of \cite{almost} 4.5). And
\[I\mathcal{O}_v\subset \mathcal{O}_v(-v)\subset \mathcal{O}_v(r_v v)=(\psi^{!}J)_v=(h^{!}I)_v\] 

So $X'$ presents at worst rational singularity at $x'$ ($I\mathcal{O}_v\subset (h^{!}I)_v$) and this rational singularity is of multiplicty $1$ (the inclusion $I\mathcal{O}_v\subset (h^{!}I)_v$ is strict). That is, $X'$ is regular, hence smooth over $\mathrm{Spec}(k)$, at $x'$.

\smallskip

\emph{m) Both morphisms $\Sigma'_i\to A$, $i=0, \infty$, are isomorphisms }:

\smallskip

Fix a basis $\lambda$ of $J^{-1}\Omega^d_{P_L/k}$. The invertible module $(I|V')^{-1}\Omega^d_{V'/k}$ of $d$-forms on $V'/\mathrm{Spec}(k)$ with simple poles along $Z'|V'$ has basis $\lambda_{V'}$, the unique extension of $\lambda|D'$ to $V'$. 
\smallskip

Denote by $p$ the composition $Z'|V'\hookrightarrow V'\to A$; it is of local complete intersection and is \'{e}tale at the two points of $Z'_t$.
Consider the fundamental class of $p$,
\[cl_p: p^*\Omega^{d-1}_{A/k}\to p^{!}\Omega^{d-1}_{A/k}=H^{-d+1}K_{Z'}|V'\] where $K_{Z'}=z^{'!}k\in\mathrm{Ob}\ D(Z', \mathcal{O})$, $z': Z'\to \mathrm{Spec}(k)$ the structural morphism.
\smallskip

With respect to the basis $\mathrm{Res}_{Z'|V'}(\lambda_{V'})$ of $H^{-d+1}K_{Z'}|V'$ and the volume form $\mathrm{Res}_0(\lambda)$ on $A$, one can write
\[cl_p(\mathrm{Res}_0(\lambda))=c\ \mathrm{Res}_{Z'|V'}(\lambda_{V'})\] with coefficient 
\[c\in H^0(Z'|V', \mathcal{O})\] 

\smallskip

This function $c$ takes value $1$ (resp. $-1$) on $\Sigma'_0|V'$ (resp. $\Sigma'_{\infty}|V'$), for $c|Z'_t$ is $1$ (resp. $-1$) at $0$ (resp. $\infty$). 
\smallskip

Naturally, $\Sigma'_0\cap\Sigma'_{\infty}\cap V'=\emptyset$, if $k$ is of characteristic $\neq 2$. 
\smallskip

In all characteristics, $cl_p$ an isomorphism. Indeed, 
\[Z'|V'=(\Sigma'_0+\Sigma'_{\infty})|V'\] by \emph{g}) and because $V'-D'$ is of codimension $\geq 2$ in $V'$. 
\smallskip

From that $cl_p$ is an isomorphism, it follows that $p$ is unramified, and that $p|(\Sigma'_i|V')$ is unramified thus \'{e}tale (SGA 1 I 9.11). 
\smallskip

By purity of branch locus \cite{almost} 2.4, if $\overline{\Sigma}'_i$ denotes the normalization of $\Sigma'_i$, the projection $\overline{\Sigma}'_i\to A$ is \'{e}tale ($\mathrm{codim}(X'-V', X')\geq 3$) and hence an isomorphism by Zariski's Main Theorem. So $\Sigma'_i\to A$ is an isomorphism.

\smallskip

Consider henceforth $\Sigma'_i$, $i=0, \infty$, as $A$-sections of $X'$. 

\smallskip

\smallskip

\emph{n) The morphism $X'\to A$ is equi-dimensional of relative dimension $1$ with $U$ dense in every fiber }:

\smallskip

Because it is smooth of relative dimension $1$ outside of $\Sigma'_0\cup\Sigma'_{\infty}$, \emph{h}), and restricts to finite morphisms on $\Sigma'_i$, $i=0, \infty$, \emph{m}).

\smallskip

\emph{o) Over each codimension $\leq 1$ point of $A$, the pseudo-isomorphism $\varphi'$ is everywhere defined and an isomorphism }:

\smallskip

This is evident on \emph{k}) and \emph{h}): over any codimension $\leq 1$ point of $A$, the two $A$-sections $\Sigma'_0$ and $\Sigma'_{\infty}$ cannot intersect.

\smallskip

\emph{p) The pseudo-isomorphism $\varphi'$ is everywhere defined and an isomorphism }:

\smallskip

Let $\Gamma$ (resp. $\overline{\Gamma}$) be the graph of $\varphi'$ (resp. the closed image of $\Gamma$ in $X'\times_AP_L$). Let $\Gamma^o$ consist of every point of $\overline{\Gamma}$ where the projection $p_2|\overline{\Gamma}: \overline{\Gamma}\to P_L$ is \'{e}tale. By Zariski's Main Theorem, $p_2|\Gamma^o$ is an open immersion. By \emph{o}) and by \cite{deligne} Note (1) 2 (ii), $\Gamma^o$ is a neighborhood of the $A$-sections $\Sigma'_i\times_A i$, $i=0, \infty$, so
\[(\varphi')^{-1}=p_1\circ (p_2|\Gamma^o)^{-1}\] is defined in a neighborhood of $V(J)$. At each point of $V(J)$, $(\varphi')^{-1}$ is obviously $A$-fiberwise quasi-finite in view of \emph{n}), hence quasi-finite, hence \'{e}tale by Zariski's Main Theorem. Therefore, $X'\to A$ is smooth also at the points of $\Sigma'_i$, $i=0, \infty$, all whose fibers are projective lines. Then by \emph{h}), $\Sigma'_0$ and $\Sigma'_{\infty}$ are disjoint. The assertion \emph{p}) follows.

\smallskip

Identify from now on $X'$ with $P_L$ via $\varphi'$, so that $I=J$, $Z'=V(J)$, $U=L^*$.

\smallskip

Let $D_X$ denote the duality auto-functor of $D^b_{\mathrm{Coh}}(X, \mathcal{O})$ with respect to its dualizing object $\mathcal{O}_X=H^{-d}f^{!}k$.

\smallskip

\emph{q) In $D^b_{\mathrm{Coh}}(X, \mathcal{O})$, one has
\[D_X(\mathcal{O}_{X'})=I,\ D_X(\mathcal{O}_{X'}/\mathcal{O}_X)=\mathcal{O}_Z[-1],\ D_X(\mathcal{O}_Z)=(\mathcal{O}_{X'}/\mathcal{O}_X)[-1]\]}
This first identity follows by definition of $I$, for $X'=P_L$ is smooth over $\mathrm{Spec}(k)$, in particular, Cohen-Macaulay. The other two are equivalent by biduality. One proves the second identity by applying $D_X$ to the exact sequence
\[0\to \mathcal{O}_X\to \mathcal{O}_{X'}\to \mathcal{O}_{X'}/\mathcal{O}_X\to 0\]

\emph{r) The algebraic space $Z$ is normal, Cohen-Macaulay and of dualizing module $\mathcal{O}_{X'}/\mathcal{O}_X=H^{-d+1}z^{!}k$, $z: Z\to \mathrm{Spec}(k)$ the structural morphism }:

\smallskip

Note first that $Z$ is reduced. In fact, $Z'$ is reduced, and one has the inclusion
\[\mathcal{O}_Z=\mathcal{O}_X/I\subset\mathcal{O}_{X'}/I=\mathcal{O}_{Z'}\]

If $\widetilde{Z}$ denotes the normalization of $Z$, then $\mathcal{O}_{\widetilde{Z}}/\mathcal{O}_Z$, which is a sub-$\mathcal{O}_Z$-module of the Cohen-Macaulay module $\mathcal{O}_{Z'}/\mathcal{O}_Z=\mathcal{O}_{X'}/\mathcal{O}_X$ and has support everywhere of codimension $\geq 1$ in $Z$, is $0$. So $Z$ is normal.

\smallskip

The other two statements together and the formula of \emph{q}),
\[D_X(\mathcal{O}_Z)=(\mathcal{O}_{X'}/\mathcal{O}_X)[-1],\] are equivalent.

\smallskip

\emph{s) Either $Z$ is integral with each projection $\Sigma'_i\to Z$ biregular, or $Z$ has two irreducible components $\rho(\Sigma'_i)$ with each extension $k(i)/k(\zeta_i)$ quadratic, where $\zeta_i$ denotes the generic point of $\rho(\Sigma'_i)$, $i=0, \infty$ }:

\smallskip

This is clear on \emph{h}) and \emph{r}): the projection $Z'\to Z$ is finite, $Z$ normal, and the module $\mathcal{O}_{Z'}/\mathcal{O}_Z=\mathcal{O}_{X'}/\mathcal{O}_X$, being dualizing on $Z$, is invertible generically.

\smallskip

\emph{t) When $Z$ is integral, $X$ is $C_L$ }:

\smallskip

That is, $X$ is pinched from $P_L$ along the finite morphism
\[\Sigma'_0\cup\Sigma'_{\infty}\to Z\] where each $A=\Sigma'_i\to Z$, $i=0, \infty$, is an isomorphism. This is clear.

\smallskip

\smallskip

\emph{u) When $Z$ has two irreducible components, $X$ is pinched from $P_L$ along the finite morphisms 
\[\rho: A=\Sigma'_i\to \Sigma_i:=\rho(\Sigma'_i)\] both of generic degree $2$, where $\Sigma_i$ is normal, integral, Cohen-Macaulay and of dualizing module $K_i=\mathcal{O}_A/\mathcal{O}_{\Sigma_i}$, $i=0, \infty$ }:

\smallskip

Again, this is clear.

\smallskip

To consider the case \emph{u}) more closely, fix $i\in\{0, \infty\}$, put $\Sigma=\Sigma_i$, identify $\Sigma'_i$ with $A$, let $K=\mathcal{O}_A/\mathcal{O}_{\Sigma}$, and choose a volume form on $A$, $\mu\in H^0(A, \Omega^{d-1}_{A/k})$.

\smallskip

--- \emph{Assume, in u), $k$ of characteristic $2$, and the extension $k(A)/k(\Sigma)$ inseparable.}

\smallskip

The inclusion of function fields
\[k(A)^2\subset k(\Sigma)\subset k(A)\] corresponds, according to Weil \cite{neron} 4.4/1, to a factorization
\[F: A\stackrel{\rho}{\to}\Sigma\stackrel{\varepsilon}{\to}A'\] of the Frobenius of $A/k$, $A'$ the base change of $A$ by the Frobenius
\[F_k: \mathrm{Spec}(k)\to \mathrm{Spec}(k),\ F_k^*: a\mapsto a^2,\ a\in k\]

Note that $\varepsilon: \Sigma\to A'$ is finite. For, it is proper and, as $F$ is finite and $\rho$ surjective, quasi-finite. This finite morphism $\varepsilon$ is also flat, since $\Sigma$ is Cohen-Macaulay and $A$ regular (EGA O 17.3.5). 

\smallskip

The module $K$, as it is dualizing on $\Sigma$, is flat over $A'$. From the exact sequence
\[0\to \mathcal{O}_{\Sigma}\to\mathcal{O}_A\to K\to 0\] one finds by duality relative to $K$ that
\[\mathcal{O}_{\Sigma}\simeq\mathrm{Ker}\ \mathrm{Tr}_{\rho}: \Omega^{d-1}_{A/k}\to K\] That is, $\mathrm{Tr}_{\rho}(\mu)=0$, hence $\mathrm{Tr}_F(\mu)=C(\mu)=0$, where $C$ denotes the Cartier operation. This implies, for example, when $d=2$, that $A$ is a supersingular elliptic curve, $\Sigma=A'$, and $K$ is $B\Omega^1_{A/k}$, the module of locally exact forms, which admits a basis over $\mathcal{O}_{A'}$.

\smallskip

--- \emph{Assume, in u), $k$ of characteristic $2$, and the extension $k(A)/k(\Sigma)$ separable.}

\smallskip

Call $\sigma$ the generator of $\mathrm{Gal}(k(A)/k(\Sigma))\simeq \mathbf{Z}/2\mathbf{Z}$. By Weil \cite{neron} 4.4/1, $\sigma$ extends to an involution of $A$. Hence, $\Sigma$ is the quotient $A/\sigma$.

\smallskip

On $\mathcal{O}_A$, one has $(\sigma^*-1)^2=(\sigma^*)^2-1=0$. And, $\sigma^*\mu=\mu$, 
\[\mathcal{O}_{\Sigma} \mu=\mathrm{Ker}\ \sigma^*-1: \Omega^{d-1}_{A/k}\to \Omega^{d-1}_{A/k}\]

When $\sigma$ has fixed points, one may by a translation assume that it fixes the origin of $A$, that is, assume that $\sigma$ is a $k$-group automorphism of $A$. Define
\[A_1=\mathrm{Im}(1+\sigma),\ A_{-1}=\mathrm{Im}(1-\sigma)\] two $\sigma$-invariant sub-abelian varieties of $A$. On $A_1$ and $A_{-1}$, $\sigma$ acts as $1$ and $-1$ respectively. The sum
\[A_1\times_SA_{-1}\to A, \ a, b\mapsto a+b\] is an isogeny and $\sigma$-compatible, which by quotient gives
\[A_1\times_S(A_{-1}/\sigma)\to A/\sigma=\Sigma\]
 
Observe that the endomorphism 
$\mathrm{1-\sigma}: A\to A$ has kernel of (pure) codimension $1$ in $A$, so that $A_{-1}$ is an elliptic curve, and $A_{-1}/\sigma$ a projective line. 
\smallskip

In fact, $\rho: A\to \Sigma$ is \'{e}tale if it is \'{e}tale at all codimension $\leq 1$ points of $A$. For, then $K=\mathcal{O}_{\Sigma}\mu$, the inclusion
\[K=\mathcal{O}_A/\mathcal{O}_{\Sigma}\simeq \mathrm{Im}(\sigma^*-1)\subset \mathrm{Ker}(\sigma^*-1)=\mathcal{O}_{\Sigma}\] is bijective, the exact sequence 
\[0\to \mathcal{O}_{\Sigma}\to\mathcal{O}_A\to K\to 0\] shows that $\mathcal{O}_A$ is locally free over $\mathcal{O}_{\Sigma}$ of rank $2$, hence $A$ faithfully flat over $\Sigma$ and $\Sigma/k$ smooth, and $A$ is \'{e}tale over $\Sigma$ by purity of branch locus.

\smallskip

--- \emph{Assume, in u), $k$ of characteristic $\neq 2$.}

\smallskip

Call $\sigma$ the generator of $\mathrm{Gal}(k(A)/k(\Sigma))\simeq\mathbf{Z}/2\mathbf{Z}$. By Weil \cite{almost} 4.4/1, $\sigma$ extends to an involution of $A$, and thus $\Sigma$ is the quotient $A/\sigma$.

\smallskip

We prove now that $\sigma^*\mu=-\mu$ :

\smallskip

\emph{Case where $\sigma$ acts freely on $A$ }:

\smallskip

Then, $A\to\Sigma$ is \'{e}tale, and $\Omega^{d-1}_{\Sigma/k}$ is the kernel of 
\[\sigma^*-1: \Omega^{d-1}_{A/k}\to \Omega^{d-1}_{A/k}\]

As $k$ is of characteristic $\neq 2$, the exact sequence
\[0\to \mathcal{O}_{\Sigma}\to \mathcal{O}_A\to \mathcal{O}_A/\mathcal{O}_{\Sigma}\to 0\] splits. In particular, $H^0(\Sigma, K)=H^0(\Sigma, \mathcal{O}_A/\mathcal{O}_{\Sigma})=0.$

\smallskip

If $\sigma^*\mu$ were $\mu$ rather than $-\mu$, then $K$, as it is isomorphic to $\Omega^{d-1}_{\Sigma/k}$, would admit a basis. This proves the claim. 

\smallskip

\smallskip

\emph{Case where $\sigma$ acts on $A$ with fixed points }:

\smallskip

One may by a translation assume that $\sigma$ is a $k$-group automorphism of $A$. Consider the sub-abelian varieties of $A$,
\[A_1=\mathrm{Im}(1+\sigma),\ A_{-1}=\mathrm{Im}(1-\sigma)\] On $A_1$ and $A_{-1}$, $\sigma$ acts as $1$ and $-1$ respectively. The sum
\[A_1\times_SA_{-1}\to A, \ a, b\mapsto a+b\] is an \'{e}tale isogeny and $\sigma$-equivariant. Let $d_{-1}=\mathrm{dim}\ A_{-1}$. The set of points of $A$ where $A\to\Sigma$ is not \'{e}tale is purely of codimension $d_{-1}$ in $A$.

\smallskip

Notice that the claim $\sigma^*\mu=-\mu$ is equivalent to that $A_{-1}$ is of \emph{odd} dimension.

\smallskip

Assume on the contrary that $d_{-1}$ is even, in particular, $\geq 2$.
Then $K$ is identical to the kernel of
\[\sigma^*-1: \Omega^{d-1}_{A/k}\to \Omega^{d-1}_{A/k}\] 

Since $k$ is of characteristic $\neq 2$, the  exact sequence
\[0\to\mathcal{O}_{\Sigma}\to \mathcal{O}_A\to \mathcal{O}_{A}/\mathcal{O}_{\Sigma}\to 0\] splits, by which $H^0(\Sigma, K)=0$. 
If $\sigma^*\mu$ were $\mu$ rather than $-\mu$, then $K$ would have a non-zero global section, a contradiction.

\smallskip

\smallskip

II) \emph{Assume $n\geq 1$ and that $S$ is the spectrum of an algebraically closed field $k$.}

\smallskip

As $X$ is pseudo-isomorphic to $L^*\times\mathbf{Z}/n\mathbf{Z}$, $X$ is reduced. Moreover, $X$ is Cohen-Macaulay.
\smallskip

Let $U$ consist of every point of $X$ where $X$ is normal. Let $\rho: X'\to X$ denote the normalization of $X$, and let $U$ be identified via $\rho$ with $\rho^{-1}(U)$. Write $I=H^0\rho^{!}\mathcal{O}$ for the conductor of $X'/X$, coherent ideal in $X'$ and in $X$, by which $Z'=X'-U'$ (resp. $Z=X-U$) is given a closed sub-algebraic space structure in $X'$ (resp. $X$).

\smallskip

As in I), one argues, on each connected component $X'_j$ of $X'$, that the pseudo-isomorphism $\varphi$ induces an $k$-isomorphism $\varphi'_j$ of $X'_j$ with $P_L$, that $Z'\cap X'_j$ is formed of two disjoint $A$-sections $\Sigma'_{j0}$ and $\Sigma'_{j\infty}$, $j\in\mathbf{Z}/n\mathbf{Z}$, and that $Z$ is normal, Cohen-Macaulay, of dualizing module $\mathcal{O}_{Z'}/\mathcal{O}_Z$. 

\smallskip

Also as in I), above each connected component $\Sigma$ of $Z$, either $\rho^{-1}(\Sigma)$ consists of two connected components $\Sigma'_{j\alpha}$, $\Sigma'_{l\beta}$, both biregular to $\Sigma$, $j, l\in\mathbf{Z}/n\mathbf{Z}$, $\alpha, \beta\in\{0, \infty\}$, or, $\rho^{-1}(\Sigma)=\Sigma'_{j\alpha}$ is connected generically of degree $2$ over $\Sigma$, $j\in\mathbf{Z}/n\mathbf{Z}$, $\alpha\in\{0, \infty\}$. 
\smallskip

These $\Sigma'_{j\alpha}$ are $2n$ in number, and $X'$ has $n$ connected components. So $X$, as it is \emph{connected} by hypothesis, is pinched from $X'$ by gluing along the $\Sigma'_{j\alpha}$,

\smallskip

--- \emph{either in a circle in which case $X=C_{nL}$},

\smallskip

--- \emph{or in a chain on whose two ends the restrictions $\rho|\Sigma'_{10}$, $\rho|\Sigma'_{n\infty}$ are generically of degree $2$.}

\smallskip

\smallskip

\smallskip

III) \emph{General case.}

\smallskip

Let $U$ consist of all points of $X$ where $f$ is smooth. Its fiber over each geometric point $s$ of $S$ is formed exactly of all points of $f^{-1}(s)$ where $f\times_Ss$ is smooth, since by assumption $f$ is flat of finite presentation.

\smallskip

Let $D=\mathrm{Dom}(\varphi)$.

\smallskip

\emph{a) The $S$-algebraic space $U$ is parafactorial along $U-D$ universally. That is, for any $S$-algebraic space $S_1$, if $U_1$ (resp. $D_1$) denotes the base change of $U$ (resp. $D$) by $S_1\to S$, $U_1$ is parafactorial along $U_1-D_1$ }:

\smallskip

For, the complement of $D\cap U$ in $U$ is fiberwise of codimension $\geq 2$ in $U$ by hypothesis of minimal model, and $U$ smooth over $S$.

\smallskip

\emph{b) The composition of the pseudo-isomorphism $\varphi$ with the projection $L^*\times\mathbf{Z}/n\mathbf{Z}\to A$ is defined on all of $U$ }:

\smallskip

In terms of the dual abelian algebraic space $A^*$ of $A$, the composition
\[D|U\stackrel{\varphi}{\to} L^*\times\mathbf{Z}/n\mathbf{Z}\to A\] can be interpreted as the data of an invertible module on the abelian algebraic space over $D|U$,
$A^*\times_S D|U$, which is fiberwise numerically trivial and trivial along the zero section. 
\smallskip

Call this invertible module $N$. By \emph{a}), an extension of $N$ to an invertible module on $A^*\times_SU$ exists, which is up to unique isomorphisms unique, trivial along the zero section, and, since $A$ is open and closed in $\mathrm{Pic}_{A^*/S}$, fiberwise numerically trivial, hence the claim.

\smallskip

\smallskip

\emph{c) The pseudo-isomorphism $\varphi$ is defined on all of $U$ }:

\smallskip

In view of \emph{b}), it suffices to observe that $L^*$ is \emph{affine} over $A$, and that, by \emph{a}), $\mathrm{prof}_{U-D}(U)\geq 2$.

\smallskip

\smallskip

\emph{d) The morphism $\varphi|U$ is an isomorphism of $U$ with $L^*\times\mathbf{Z}/n\mathbf{Z}$ }:

\smallskip

The assertion can be verified fiber by fiber. Then it is immediate on the study made in I) and II).  

\smallskip

\emph{e) The geometric number $N(s)$ of connected components of $f^{-1}(s)-\mathrm{Im}(\varphi^{-1}_s)$, $s\in S$, is a constructible function on $S$ }:

\smallskip

This follows from EGA IV 9.7.9.

\smallskip

\emph{f) The function $s\mapsto N(s)$ is upper semi-continuous on $S[1/2]$ }:

\smallskip

By \emph{e}), one may assume that $S$ is the spectrum of a discrete valuation ring with generic point $t$ and closed point $s$, $k(s)$ algebraically closed. One needs show that the case where $N(t)=n+1$, $N(s)=n$ does not happen. 

\smallskip

Otherwise, replacing if necessary $S$ by a spectrum $S'$ of discrete valuation ring faithfully flat over $S$, and $(f, \varphi)$ by $(f\times_SS', \varphi\times_SS')$, one would have $f^{-1}(s)=C_{nL}$, and that $f^{-1}(t)-\mathrm{Im}(\varphi^{-1}_t)$ admit one component of the form $A_t/\sigma_t$, where $\sigma_t$ is an involutive automorphism of $A_t$. Recall from I) that $\sigma^*_t$ acts as $-1$ on the volume forms of $A_t$.

\smallskip

As $A$ is the $S$-N\'{e}ron model of $A_t$, a unique involution $\sigma$ of $A$ extends $\sigma_t$. Write
\[\sigma(x)=\alpha(x)+a,\ x\in A\] where $\alpha$ is an $S$-group automorphism of $A$, and $a\in A(S)$. That $\sigma^2=1$ amounts to that $\alpha^2=1$, $\alpha(a)=-a$. 
\smallskip

Let $\Sigma$ denote the closed image of $A_t/\sigma_t$ in $X$. Its closed fiber $\Sigma_s$ is connected by Zariski and purely of dimension $d-1$, hence is a connected component of $f^{-1}(s)-\mathrm{Im}(\varphi^{-1}_s)$. In particular, $(\Sigma_s)_{\mathrm{red}}\simeq A_s$ and does not contain $k(s)$-rational curves, so the projection $A_t\to \Sigma_t$ extends to an $S$-morphism $A\to \Sigma$ (\cite{almost} 2.1).

\smallskip

This morphism $A\to \Sigma$ is proper, thus surjective, thus finite. For, $A_s\to (\Sigma_s)_{\mathrm{red}}=A_s$, being surjective, is up to translation an isogeny, with say finite kernel $E_s$. 

\smallskip

Therefore, $\Sigma$ is of normalization $A/\sigma$. Note that $(A/\sigma)_s=A_s/\sigma_s$, since $k(s)$ is of characteristic $\neq 2$. The factorization 
\[A_s\to A_s/\sigma_s\to (\Sigma_s)_{\mathrm{red}}=A_s\] implies that $\sigma(x)-x$, $\forall\ x\in A_s(s)$, are contained in $E_s(s)$.
Namely, the sub-abelian variety $\mathrm{Im}(\alpha_s-1)$ of $A_s$ is contained in $-a_s+E_s$. 
\smallskip

But $E_s$ is finite over $s$. So $\mathrm{Im}(\alpha_s-1)=0$, $\alpha_s=1$, and so $\alpha=1$, since the functor of $S$-group automorphisms of $A$ is representable and locally unramified over $S$.
\smallskip

This shows that $\sigma$ is a translation. Hence, $\sigma^*_t\mu$ is $\mu$ rather than $-\mu$ for any volume form $\mu$ on $A_t$, a contradiction.

\smallskip

\smallskip

\emph{g) The function $s\mapsto N(s)$ is lower semi-continuous on $S[1/2]$ }:

\smallskip

Since $s\mapsto h^1(s):=\mathrm{dim}_{k(s)}\ H^1(f^{-1}(s), \mathcal{O})$ is upper semi-continuous, it suffices to show that $h^1(s)$ is $d(s)$ and $<d(s)$, when $N(s)=n$ and $N(s)=n+1$ respectively, where $d(s)=\mathrm{dim}\ f^{-1}(s)$, $s\in S[1/2]$; note that $f$ is equi-dimensional, as $\mathrm{Im}(\varphi^{-1})$ is $S$-dense in $X$.

\smallskip

--- \emph{Calculation of $h^i(s):=\mathrm{dim}_{k(s)}\ H^i(f^{-1}(s), \mathcal{O})$, when $N(s)=n$ }:

\smallskip

One may assume $S=s=\mathrm{Spec}(k)$, $k$ an algebraically closed field. Then $X=C_{nL}$. Let $d=\mathrm{dim}\ X$. Consider the canonical exact sequence
\[0\to \mathcal{O}_X\to \mathcal{O}_{X'}\to \mathcal{O}_{X'}/\mathcal{O}_X\to 0\] where $X'=P_L\times \mathbf{Z}/n\mathbf{Z}$ is the normalization of $X$. 
\smallskip

With the notations of II), write 
\[\Pi: \mathcal{O}_{X'}\to \mathcal{O}_{X'}/\mathcal{O}_X=\mathcal{O}_{Z'}/\mathcal{O}_Z\] for the projection.

\smallskip

Both $H^0(X', \mathcal{O}_{X'})$ and $H^0(X, \mathcal{O}_{X'}/\mathcal{O}_X)=H^0(Z, \mathcal{O}_{Z'}/\mathcal{O}_Z)$ are of dimension $n$ over $k$, and the exact sequence
\[0\to H^0(\mathcal{O}_X)\to H^0(\mathcal{O}_{X'})\to H^0(\mathcal{O}_{X'}/\mathcal{O}_X)\] shows that $\mathrm{Ker}\ H^0(\Pi)$ is of dimension $1$. So $\mathrm{Coker}\ H^0(\Pi)$ is of dimension $1$ as well.

\smallskip

Observe that, for each integer $i$, $H^i(A, \mathcal{O})\simeq H^i(P_L, \mathcal{O})$, via which the map
\[H^i(\Pi): H^i(X', \mathcal{O}_{X'})\to H^i(X, \mathcal{O}_{X'}/\mathcal{O}_X)\] can be identified with
\[H^i(A, \mathcal{O})\otimes H^0(\Pi)\] In particular, both kernel and cokernel of $H^i(\Pi)$ have the same dimension as $H^i(A, \mathcal{O})$. From the cohomology exact sequence
\[0\to \mathrm{Coker}\ H^{i-1}(\Pi)\to H^i(X, \mathcal{O})\to \mathrm{Ker}\ H^i(\Pi)\to 0\] it follows finally that $H^i(X, \mathcal{O})$ is of dimension
\[\binom{d-1}{i-1}+\binom{d-1}{i}=\binom{d}{i}\]

--- \emph{Proof of $h^1(s)<d(s)$, when $N(s)=n+1$ }:

\smallskip

Similarly, one may assume $S=s=\mathrm{Spec}(k)$, $k$ an algebraically closed field of characteristic $\neq 2$. Let $d=\mathrm{dim}\ X$. Again, consider the exact sequence
\[0\to \mathcal{O}_X\to \mathcal{O}_{X'}\to \mathcal{O}_{X'}/\mathcal{O}_X\to 0\]

In the notations of II), $\mathcal{O}_{X'}/\mathcal{O}_X=\mathcal{O}_{Z'}/\mathcal{O}_Z$.

\smallskip

Recall that $X$ is pinched from $X'=P_L\times\mathbf{Z}/n\mathbf{Z}$ in a chain. On its two ends, $\rho(\Sigma'_{10})=\Sigma_{10}$, $\rho(\Sigma'_{n\infty})=\Sigma_{n\infty}$, the dualizing module $K_Z$ has no non-zero global sections; in between, which consists of $n-1$ connected components, $K_Z$ is isomorphic to $\mathcal{O}_Z$. So $H^0(Z, \mathcal{O}_{Z'}/\mathcal{O}_Z)=H^0(Z, K_Z)$ is of dimension $n-1$ over $k$, thus
\[0\to H^0(X, \mathcal{O}_{X})\to H^0(X', \mathcal{O}_{X'})\to H^0(Z, \mathcal{O}_{Z'}/\mathcal{O}_Z)\to 0\] is exact, and thus 
\[H^1(X, \mathcal{O})\simeq \mathrm{Ker}\ H^1(X', \mathcal{O})\to H^1(Z, \mathcal{O}_{Z'}/\mathcal{O}_Z)\] which as one verifies easily is isomorphic to
\[\mathrm{Ker}\ H^1(X'_1, \mathcal{O})\to H^1(\Sigma_{10}, \mathcal{O}_{Z'}/\mathcal{O}_Z)\] where $X'_1=P_L\times \{1\}$ is the left most component of $X'$. This finishes the proof, in view that $H^1(X'_1, \mathcal{O})\simeq H^1(A, \mathcal{O})$ is of dimension $d-1$ over $k$.

\smallskip

\emph{h) The function $s\mapsto N(s)$ is locally constant on $S[1/2]$ }:

\smallskip

Combine \emph{f}) and \emph{g}). 
\smallskip

The proof is now complete.


\bibliographystyle{amsplain}


\end{document}